\documentclass{article}

\usepackage[ansinew]{inputenc}
\usepackage{amsmath}
\usepackage{amssymb}
\usepackage{amsthm}
\usepackage{amscd}
\usepackage{amsfonts}
\usepackage[english]{babel}
\usepackage{graphicx}
\usepackage{cancel}

\newtheorem {thm}{Theorem}[section]

\newtheorem {defn}[thm]{Definition}

\newtheorem {lemma}[thm]{Lemma}

\newtheorem {cor}[thm]{Corollary}

\theoremstyle{remark}
\newtheorem{ExInternal}{Exercise}[section]

\makeatletter
\let\@exercises\@empty%
\newcommand\exercise[2][]{%
    \g@addto@macro\@exercises{%
        \begin{ExInternal}[#1]%
            #2%
        \end{ExInternal}%
    }%
}

\newcommand\exerciseshere{%
    \subsection*{Exercises}
    \@exercises%
    \global\let\@exercises\@empty%
}
\makeatother

\newcommand{\I}{\boldsymbol{\mathcal{I}}}

\newcommand{\SU}{\boldsymbol{\mathcal{SU}}}
\newcommand{\U}{\boldsymbol{\mathcal{U}}}

\newcommand{\R}{\boldsymbol{\mathcal{R}}}

\newcommand{\TR}{\operatorname{tr}}

\newcommand{\bbmatrix}[1]{\left[ \begin{array}{ccccccccccccccccccc} #1 \end{array} \right]}
\newcommand{\SPEC}{{\operatorname{spec}}}

\newcommand{\Rems}{\noindent\textbf{Remarks: 1. }}
\newcommand{\Proof}{\noindent\textbf{Proof: }}

\title{Inverse scattering of Canonical systems and their evolution}
\author{Andrey Melnikov\\Drexel University}

\begin{document}
\maketitle
\abstract{In this work we present an analogue of the inverse scattering for Canonical systems using theory of vessels and associated to them
completely integrable systems. Analytic coefficients fits into this setting, significantly expanding the class of functions for which the inverse
scattering exists. We also derive an evolutionary equation, arising from canonical systems, which describes an evolution of
the logarithmic derivative of the tau function, associated to these systems.}
\tableofcontents

\section{Introduction}
Starting from the works of M. Krein \cite{bib:KreinCanSys, bib:KreinMelikCanSys}, Gasymov \cite{bib:GasymobCanSys} and continuing with research of Alpay-Gohberg \cite{bib:agCanSys}, 
Dym-Arov \cite{bib:ArDymCanSys} and many other mathematicians, the canonical systems
where one of the basic models, studied in order to understand completely integrable Partial Differential Equations (PDEs). A system of equations
\[ \bbmatrix{0&-1\\1&0} \dfrac{\partial}{\partial x} y = \mu y+ \bbmatrix{\widetilde p & \widetilde q \\\widetilde q &-\widetilde p} y, \quad 
y=\bbmatrix{y_1(\mu,x)\\y_2(\mu,x)}, \mu\in\mathbb C
\]
for the vector $y$, is called the \textit{Canonical System} and $\mu$ is called the \textit{spectral parameter}. This system and its generalizations were
studied for various classes of functions: Wienner class \cite{bib:ArDymCanSys}, rational functions \cite{bib:agCanSys}, rational spectral function 
\cite{bib:GKSCanSys}, strictly preudo-exponential \cite{bib:AGKSCanSys}. It is worth noticing that most of the techniques to solve and research Canonical Systems
are related to the theory of systems, $J$-contractive functions, reproducing kernel Hilbert spaces and related topics.

We are going to present a new, related to these techniques, approach using a theory of vessels. The main advantage of this approach is that we are able to present the 
inverse scattering for arbitrary analytic coefficients $\widetilde p, \widetilde q$. Moreover, a natural notion of KdV vessels enables to use these ideas
to solve an initial value problem for the following completely integrable PDE, which we will call \textit{canonical PDE}
\begin{equation} \label{eq:CanSysEvol}
 \dfrac{\partial^2}{\partial t^2} \beta = \dfrac{\partial}{\partial x} [ - \dfrac{1}{2} (\beta')^2 - \dfrac{1}{4}\beta''' + \dfrac{(\dfrac{\partial}{\partial t} \beta)^2+ \dfrac{1}{4}(\beta'')^2}{\beta'}].
\end{equation}
The idea behind this solution is the usual inverse scattering method, developed in \cite{bib:GGKM} for the KdV equation. Since we are able to define the scattering
data on $\mathbb R$ for analytic parameters, we are also able to solve this initial value problem. The name of KdV vessel, which serves as a solution to this initial value PDE
comes from the analogous problem for the KdV equation. We use the same ideas and constructions. Still, the derivation of the new formulas, based on the same principal
construction of a vessel requires some research, which is presented here and describes a solution of a very important object, Canonical Systems.

The setting, developed for vessels in \cite{bib:amv, bib:SchurIEOT, bib:SLVessels, bib:MelVinC} requires slightly different operators, then in the original definition of Canonical systems.
if we multiply the original equation by $i$, we can rewrite it as
\begin{equation} \label{eq:CanSys}
\lambda y - \bbmatrix{0&i\\-i&0} \dfrac{\partial}{\partial x} y + \bbmatrix{i p & iq \\iq &-ip} y = 0.
\end{equation}
and, as a result, the scattering of the operator
\begin{equation} \label{eq:DefL}
L = - \bbmatrix{0&i\\-i&0}\dfrac{\partial}{\partial x} + \bbmatrix{i p & iq \\iq &-ip}
\end{equation}
will be our primary interest. We are going to consider a scattering theory in a broader sense, since we will show that solutions of
\eqref{eq:CanSys} with trivial potential ($p=q=0$) will be mapped to solutions of $\eqref{eq:CanSys}$ with deriived from the vessel 
parameters $p(x),q(x)$, except for values of $\lambda$,
lying in a spectrum of a certain operator. An analogous scattering scheme is presented in Section \ref{sec:ClassicScatt} and coincides
with the classical inverse scattering when the potential
$\gamma_*(x)$ is classical one (namely satisfies $\int\limits_{\mathbb R} |p(x)|dx, \int\limits_{\mathbb R} |q(x)|dx < \infty$).
Main contribution of this paper is that also analytic initial potentials (when $p(x)$, $q(x)$
are analytic functions) can serve as initial values for \eqref{eq:CanSys} and there is a prescribed way how to solve it, using vessels.
In order to do it, a more general setting of evolutionary
nodes (prevessels) is presented. We show in Section \ref{sec:AnalScatt} how the formulas are incorporated to create this. Moreover,
when served as an initial condition, we show that a solution of \eqref{eq:CanSysEvol} is created, which coincides for $t=0$ with the
given value, which is created from the scattering of $L$ \eqref{eq:DefL}.

The notion of KdV vessel comes from an analogous solution of the Korteweg-de Vries (KdV) and evolutionary Non Linear Schr\" odinger (NLS) equations. Consider the following collection of operators and spaces
\[ \begin{array}{lllllllll}
\mathfrak{V}_{KdV} = \bbmatrix{C(x) & A_\zeta, \mathbb X(x), A & B(x)&  \sigma_1,\sigma_2,\gamma,\gamma_*(x) \\
\mathbb C^2 & \mathcal{K} & \mathbb C^2 & \Omega}, \\
\hspace{2cm} B, C^*:\mathbb C^2\rightarrow\mathcal K, \quad A,\mathbb X,A_\zeta:\mathcal K\rightarrow\mathcal K, \\
\hspace{2cm} \sigma_1, \sigma_2, \gamma,\gamma_*(x, t) \text{ - $2\times 2$ matrices}, \quad \Omega\subseteq\mathbb R^2\\
\dfrac{\partial}{\partial x} B \sigma_1 = - (A B \sigma_2 + B \gamma), \hspace{1.2cm} \dfrac{\partial}{\partial t} B = i A  \dfrac{\partial}{\partial x} B, \\
\dfrac{\partial}{\partial x} C^* \sigma_1  = - (A_\zeta^* C^* \sigma_2 + C^* \gamma), \quad \quad \dfrac{\partial}{\partial t} C^*  = i A_\zeta^*  \dfrac{\partial}{\partial x} C^*, \\
\dfrac{\partial}{\partial x} \mathbb X =  B \sigma_2 C, \hspace{3cm} \dfrac{\partial}{\partial t} \mathbb X = i A B \sigma_2 C - i B \sigma_2 C A_\zeta + i B \gamma C, \\ 
A \mathbb X + \mathbb X A_\zeta  + B \sigma_1 C  = 0, \\
\gamma_* =  \gamma + \sigma_2 C \mathbb X^{-1} B \sigma_1 - \sigma_1 C \mathbb X^{-1} B \sigma_2.
\end{array} \]
which is called a \textbf{KdV vessel} on $\Omega$, where $\mathbb X(x, t)$ is invertible. We obtain a solution of the KdV equation \eqref{eq:KdV} on $\Omega$ as follows. For 
\[ \sigma_1 = \bbmatrix{0&1\\1&0}, \quad \sigma_2=\bbmatrix{1&0\\0&0}, \quad \gamma=\bbmatrix{0&0\\0&i}
\]
the function ($(x_0,t_0)\in\Omega$)
\[ q_{\mathfrak V}(x,t) = -2\dfrac{\partial^2}{\partial x^2} \ln \det(\mathbb X^{-1}(x_0,t_0)\mathbb X(x,t))
\]
satisfies \cite{bib:KdVVessels, bib:AnalPotKdV}
\begin{equation} \label{eq:KdV}
q_t = -\dfrac{3}{2} q q_x + \dfrac{1}{4} q_{xxx}.
\end{equation}
If we choose vessel parameters 
\[ \sigma_1 = I= \bbmatrix{1&0\\0&1}, \quad \sigma_1 = \dfrac{1}{2}\bbmatrix{1&0\\0&-1}, \quad \gamma = 0 = \bbmatrix{0&0\\0&0},
\]
then there created solutions $y=\bbmatrix{0&1}\gamma_*(x,t)\bbmatrix{1\\0}$ of the evolutionary Non Linear Schr\" odinger (NLS) equation
\begin{equation} \label{eq:ENLS}
i y_t + y_{xx} + 2 |y|^2 y = 0.
\end{equation}
The proof of this fact can be found in \cite{bib:ENLS}.

In this paper, as we already pointed out, we study the same scheme (the same vessels) to study solutions of the
canonical PDE \eqref{eq:CanSysEvol}. For this it is necessary to choose 
\[ \sigma_1 = \bbmatrix{0 & i \\ -i &  0},
   \sigma_2 = I = \bbmatrix{1 & 0 \\ 0 & 1},
   \gamma  = 0 = \bbmatrix{0 & 0 \\ 0 & 0}.
\]
Then the function $\beta(x,t) = \dfrac{\partial}{\partial x} \tau(x,t) = \dfrac{\partial}{\partial x} \det (\mathbb X^{-1}(0,0)\mathbb X(x,t))$
satisfies \eqref{eq:CanSysEvol}.

% -------------------------------------------------------------------------------------------
\section{Background}
\subsection{Node, prevessel, vessel}
\begin{defn} A \textbf{node} is a collection of operators and spaces
\[ \mathfrak N = \bbmatrix{
C & A_\zeta, \mathbb X, A & B & \sigma_1\\
\mathbb C^2 & \mathcal K & \mathbb C^2}
\]
where $\mathcal K$ is a Krein space, $C:\mathcal K\rightarrow\mathbb C^2$,
$\mathbb X: \mathcal K\rightarrow\mathcal K$,  $B:\mathbb C^2\rightarrow\mathcal K$ are bounded operators,
$\sigma_1=\sigma_1^*$ - invertible $2\times 2$ matrix, 
$A, A_\zeta$ are generators of $C_0$ groups on $\mathcal K$ with identical dense domain $D(A)=D(A_\zeta)$. The operator $\mathbb X$
is assumed to satisfy $\mathbb X(D(A))\subseteq D(A)$. The operators of the node are subject to the \textbf{Lyapunov equation}
\begin{equation} \label{eq:Lyapunov}
A \mathbb X u + \mathbb X A_\zeta u + B \sigma_1 C u = 0, \quad \forall u\in D(A_\zeta)=D(A).
\end{equation}
If $\mathbb X$ is invertible, the \textbf{transfer function} of $\mathfrak N$ is
\begin{equation} \label{eq:S0realized}
S(\lambda) = I - C \mathbb X^{-1} (\lambda I - A)^{-1} B \sigma_1.
\end{equation}
The node $\mathfrak N$ is called \textbf{symmetric} if $A_\zeta=A^*$ and $C=B^*$.
\end{defn}
\Rems function $S(\lambda)$, representable in the form \eqref{eq:S0realized} is called \textit{realized} \cite{bib:bgr}.
\textbf{2.} if $A_\zeta=A+T$ for a bounded operator $T$, then $D(A)=D(A_\zeta)$ ($A_\zeta$ is called a perturbation of $A$ in this case).
\textbf{3.} if $\mathbb X=I$, then the condition $\mathbb X(D(A))\subseteq D(A)$ holds. \textbf{4.} when the node is symmetric one can verify that
\[ S(\lambda) \sigma_1^{-1} S^*(-\bar\lambda) = \sigma_1^{-1}
\]
at all points of analyticity of $S$. \textbf{5.} for the unbounded operators $A,A_\zeta$ to be generators of $C_0$-groups, it necessary and sufficient
to demand that they satisfy the conditions of the Hille-Yosida Theorem \cite[Theorem 3.4.1]{bib:Staffans}, \cite[Theorem II.3.5]{bib:EngelNagel}.
Particularly, they must be closed, densely defined operators.

In the case $\mathbb X$ is invertible, we consider a stronger notion of a node as follows.
\begin{defn} A node $\mathfrak N$ is called \textbf{invertible}, if $\mathbb X$ is invertible and $\mathbb X^{-1}(D(A))\subseteq D(A)$.
\end{defn}
A simple chain of inclusions for an invertible node
\[ \mathbb X(D(A))\subseteq D(A) \Rightarrow D(A)\subseteq \mathbb X^{-1}(D(A)) \subseteq D(A),
\]
where the first inclusion comes from the node condition, and the last one from the invertible node condition, implies that
$\mathbb X^{-1}(D(A)) = D(A)$. Similarly, $\mathbb X(D(A)) = D(A)$. Moreover, taking $u=\mathbb X^{-1}u'$, where $u,u'\in D(A)$ and plugging it into the
Lyapunov equation \eqref{eq:Lyapunov}, we obtain that
\begin{equation} \label{eq:Lyapunov-1}
A_\zeta \mathbb X^{-1} u' + \mathbb X^{-1} A u' + \mathbb X^{-1}B \sigma_1 C \mathbb X^{-1}u' = 0, \quad \forall u'\in D(A),
\end{equation}
after multiplying by $\mathbb X^{-1}$ from the left. From the existence of this Lyapunov equation we obtain the following Lemma. 
\begin{lemma} If $\mathfrak N$ is an invertible node, then
\[ \mathfrak N_{-1} = \bbmatrix{
C\mathbb X^{-1} &  A, \mathbb X^{-1},A_\zeta & \mathbb X^{-1}B & \sigma_1\\
\mathbb C^2 & \mathcal K & \mathbb C^2}
\]
is also a node.
\end{lemma}
One could consider a similar notion of ``adjointable'' node, for which the adjoint of the Lyapunov equation \eqref{eq:Lyapunov} would define
a node, but we do not insert details here. Actually, there is a theory of construction of new such nodes from old ones, similarly to the theory
presented in \cite{bib:Brodskii, bib:bgr}. In the case $\mathbb X = I$ (the identity operator) we have a very well developed theory \cite{bib:BL} of 
(symmetric) nodes with $A_\zeta = A^*$, which has a finite dimensional real part: $A + A^* = -B\sigma_1B^*$. 

Finally, rewriting the Lyapunov equation \eqref{eq:Lyapunov-1}, of the invertible node as follows
\[ (-A_\zeta) \mathbb X^{-1} u' + \mathbb X^{-1} (-A) u' + \mathbb X^{-1}B \sigma_1^{-1} \sigma_1(-\sigma_1 C) \mathbb X^{-1}u' = 0, \quad \forall u'\in D(A)
\]
we arrive to the node
\[ \mathfrak N^{-1} = \bbmatrix{
-\sigma_1 C\mathbb X^{-1} &  -A, \mathbb X^{-1}, -A_\zeta &  \mathbb X^{-1}B\sigma_1^{-1} & \sigma_1\\
\mathbb C^2 & \mathcal K & \mathbb C^2}
\]
whose transfer unction
\[ S^{-1}(\lambda) = I + \sigma_1 C (\lambda I + A_\zeta)^{-1}\mathbb X^{-1}B
\]
is the inverse of the transfer function, defined in \eqref{eq:S0realized},
of the original invertible node $\mathfrak N$. This is a standard fact, related to Schur complements and can be found in \cite{bib:Brodskii, bib:bgr}.
\begin{defn} Class $\R(\sigma_1)$ consist of $2\times 2$ matrix-valued functions $S(\lambda)$ of the complex variable $\lambda$, which are transfer functions of invertible nodes.
The subclass $\U(\sigma_1)\subseteq\R(\sigma_1)$ consists of the transfer functions of symmetric, invertible nodes.
The Schur class $\SU(\sigma_1)\subseteq\U$ demand also that the inner space $\mathcal K$ is Hilbert and $\mathbb X>0$.
\end{defn}
When $S(\lambda)$ is just analytic at infinity (hence $A$ must be bounded), 
there is a very well known theory of realizations developed in \cite{bib:bgr}. 
For analytic at infinity and symmetric, i.e. satisfying $S^*(-\bar\lambda) \sigma_1 S(\lambda) =  \sigma_1$,
functions there exists a good realization theory using Krein spaces,
developed in \cite{bib:KreinReal}\footnote{At the paper \cite{bib:KreinReal} a similar result is proved for functions symmetric with respect to the unit circle,
but it can be translated using Calley transform into $S^*(-\bar\lambda) \sigma_1 S(\lambda) =  \sigma_1$ and was done in \cite{bib:GenVessel, bib:SchurIEOT}}.
Such a realization is then translated into a function in $\U(\sigma_1)$.
The sub-classes $\U, \SU$ appear a lot in the literature and correspond to the symmetric case. We will not particularly
consider these two classes here and refer to \cite{bib:SchurIEOT}.

Equations, which arise in the theory of vessels involve differential equations
with unbounded operators. As a result, an operator satisfying such an equation must satisfy a relation
with the domain of the unbounded operator, which is presented in the next Definition.
\begin{defn} A bounded operator $B:\mathbb C^2\rightarrow\mathcal K$ is called $A$-\textbf{regular},
where $A: \mathcal K\rightarrow\mathcal K$ is linear, if $Be\in D(A)$ for all $e\in\mathbb C^2$.
\end{defn}
\begin{defn} \label{def:preVessel} Fix $2\times 2$ matrices $\sigma_2=\sigma_2^*$, $\gamma = -\gamma^*$.
The collection of operators and spaces 
\begin{equation} \label{eq:DefpreV}
\mathfrak{preV} = \bbmatrix{C(x) & A_\zeta, \mathbb X(x), A & B(x)&  \sigma_1,\sigma_2,\gamma \\
\mathbb C^2 & \mathcal{K} & \mathbb C^2}
\end{equation}
is called a (non-symmetric) \textbf{prevessel}, if the following conditions hold:
1. $\mathfrak{preV}$ is a node for all $x\in\mathbb R$, 2.
the operator $B(x)\sigma_2$ is $A$-regular, 3. $C(x), \mathbb X(x), B(x)$ are differentiable, bounded
operators, subject to the following conditions
\begin{eqnarray}
\label{eq:DB} \frac{\partial}{\partial x} B(x)  = - (A B \sigma_2 + B \gamma) \sigma_1^{-1}, \\
\label{eq:DC} \frac{\partial}{\partial x} C(x)u = \sigma_1^{-1}(- \sigma_2 C A_\zeta u + \gamma Cu) , \quad \forall u\in D(A_\zeta), \\
\label{eq:DX} \frac{\partial}{\partial x} \mathbb X =  B \sigma_2 C,
\end{eqnarray}
The prevessel $\mathfrak{preV}$ is called \textbf{symmetric} if
$A_\zeta=A^*$ and $C(x)=B^*(x)$ for all $x\in\mathbb R$.
\end{defn}
It turns out that the structure of a prevessel implies the Lyapunov equations \eqref{eq:Lyapunov}, \eqref{eq:Lyapunov-1} as the following
Lemma claims.
\begin{lemma}[\textbf{permanence of the Lyapunov equations}] \label{lemma:Redund}
Suppose that $B(x), C(x), \mathbb X(x)$ satisfy \eqref{eq:DB}, \eqref{eq:DC}, \eqref{eq:DX} respectively and $\mathbb X(x)(D(A_\zeta))\subseteq D(A)$
for all $x\in\mathbb R$. Then if the Lyapunov equation \eqref{eq:Lyapunov}
holds for a fixed $x_0$, then it holds for all $x$. In the case the operator $\mathbb X(x)$ is invertible and $B(x), C(x), \mathbb X(x)$ are part of an invertible node,
if \eqref{eq:Lyapunov-1} holds for a fixed $x_0$, then it holds for all $x$.
\end{lemma}
\noindent\textbf{Proof:} 
Let us differentiate the right hand side of the Lyapunov equation \eqref{eq:Lyapunov}:
\[ \begin{array}{lll}
\dfrac{d}{dx} [A \mathbb X u + \mathbb X A_\zeta u + B \sigma_1 C u ] = \\
= A B(x)\sigma_2C(x)  u + B(x)\sigma_2C(x)  A_\zeta u - A B(x) \sigma_2 C(x) u - B(x) \sigma_2 C(x) A_\zeta u \\
= 0.
\end{array} \]
The terms involving $\gamma$ are canceled, because $\gamma+\gamma^*=0$, by the assumption on it.
Thus it is a constant and the result follows. For the invertible node case, the condition \eqref{eq:Lyapunov-1} is a result of \eqref{eq:Lyapunov}. \qed

\begin{defn} \label{def:Vessel} The collection of operators, spaces and a set $\Omega\subseteq\mathbb R$ 
\begin{equation} \label{eq:DefVGen}
\mathfrak{V} = \bbmatrix{C(x) & A_\zeta, \mathbb X(x), A & B(x)&  \sigma_1,\sigma_2,\gamma,\gamma_*(x) \\
\mathbb C^2 & \mathcal{K} & \mathbb C^2 & \Omega}
\end{equation}
is called a (non-symmetric) \textbf{vessel}, if $\mathfrak{V}$ is a pre-vessel, $\mathbb X(x)$ is invertible on
$\Omega$, and $\mathfrak{V}$ is also an invertible node for all $x\in\Omega$.
The $2\times 2$ matrix-function $\gamma_*(x)$ is assumed to satisfy the \textbf{linkage condition} on $\Omega$
\begin{equation} \label{eq:Linkage}
	\gamma_*  =  \gamma + \sigma_2 C(x) \mathbb X^{-1}(x) B(x) \sigma_1 - \sigma_1 C(x) \mathbb X^{-1}(x) B(x) \sigma_2.
\end{equation}
\end{defn}

The class of the transfer functions of vessels is defined as follows
\begin{defn} Class $\I=\I(\sigma_1,\sigma_2,\gamma;\Omega)$ consist of $2\times 2$ matrix-valued (transfer) functions $S(\lambda,x)$ of the complex variable $\lambda$ and $x\in\Omega\subseteq\mathbb R$, possessing the following representation:
\begin{equation} \label{eq:DefS}
 S(\lambda,x) = I - C(x) \mathbb X^{-1}(x) (\lambda I - A)^{-1} B(x) \sigma_1,
\end{equation}
where the operators $C(x), \mathbb X(x), B(x)$ are part of a vessel $\mathfrak V$.
\end{defn}
Before we prove the B\" acklund transformation Theorem \ref{thm:Backlund} we present a technical lemma.
\begin{lemma} \label{lemma:X-1B} Let $\mathfrak{V}$ be a vessel. Then for all $u\in D(A)$
\begin{eqnarray}
\label{eq:DCX-1} \sigma_1 \dfrac{d}{dx} [C(x)\mathbb X^{-1}(x)] u = \sigma_2 C(x)\mathbb X^{-1}(x) A u + \gamma_*(x) C(x)\mathbb X^{-1}(x) u, \\
\label{eq:DX-1B} \dfrac{d}{dx} [\mathbb X^{-1}(x)B(x)] \sigma_1 = A_\zeta \mathbb X^{-1}(x) B(x)\sigma_2 - \mathbb X^{-1}(x) B(x) \gamma_*(x).
\end{eqnarray}
\end{lemma}
\Proof Consider \eqref{eq:DCX-1} first. We write under each equality the corresponding equation that is used to derive the next line:
\[ \begin{array}{lllllll}
\sigma_1 \dfrac{\partial}{\partial x} [C\mathbb X^{-1}] u & = \sigma_1 \dfrac{\partial}{\partial x} (C) \mathbb X^{-1} u +\sigma_1 C\dfrac{\partial}{\partial x} (\mathbb X^{-1}) u \\
 &  \quad \quad \text{\eqref{eq:DC}: } \quad \quad \sigma_1 \frac{\partial}{\partial x} C(x)u = - \sigma_2 C A_\zeta u + \gamma C u \\
 & \quad \quad \text{\eqref{eq:DX}:} \quad \quad \frac{\partial}{\partial x} \mathbb X =  B \sigma_2 C \\
 & = - \sigma_2 C A_\zeta \mathbb X^{-1} u + \gamma C \mathbb X^{-1} u  - \sigma_1 C\mathbb X^{-1} B \sigma_2 C \mathbb X^{-1} u \\
 & = -\sigma_2 C A_\zeta \mathbb X^{-1} u  + ( \gamma - \sigma_1 C\mathbb X^{-1} B \sigma_2) C \mathbb X^{-1} u \\
 & \quad \quad \text{\eqref{eq:Lyapunov-1}: }  \quad \quad  A_\zeta \mathbb X^{-1} u + \mathbb X^{-1} A u + \mathbb X^{-1}B \sigma_1 C \mathbb X^{-1}u = 0 \\
 & = \sigma_2 C \mathbb X^{-1} Au  + ( \gamma  + \sigma_2 C\mathbb X^{-1} B \sigma_1 - \sigma_1 C\mathbb X^{-1} B \sigma_2) C \mathbb X^{-1} u \\
 & \quad \quad \text{\eqref{eq:Linkage}: } \quad \quad  \gamma_*  =  \gamma + \sigma_2 C \mathbb X^{-1} B \sigma_1 - \sigma_1 C \mathbb X^{-1} B \sigma_2 \\
 & = \sigma_2 C\mathbb X^{-1} A u + \gamma_* C \mathbb X^{-1}u.
\end{array} \]
Notice that all equations of the vessel can be used, since we apply them to a vector $u$ from $D(A)$.

The equation \eqref{eq:DX-1B} is proved in exactly the same manner. \qed

Now we have all the ingredients of the following Theorem. This theorem has its origins at the work of M. Liv\c sic \cite{bib:Vortices} and was
proved for bounded operators in \cite{bib:GenVessel, bib:ENLS, bib:SchurIEOT}. Now we present a generalization of these results for the case of unbounded
operator $A$.
\begin{thm}[Vessel=B\" acklund transformation] \label{thm:Backlund}
Let $\mathfrak{V}$ be a vessel defined in \eqref{eq:DefVGen} and satisfying the conditions of Definition \ref{def:Vessel}. Fix
$\lambda\not\in\SPEC(A)$ and let $u(\lambda,x)$ be a solution of the input LDE 
\begin{equation} \label{eq:InCC}
 \lambda \sigma_2 u(\lambda, x) -
\sigma_1 \frac{\partial}{\partial x}u(\lambda,x) +
\gamma u(\lambda,x) = 0.
\end{equation} 
Then the function $y(\lambda,x)=S(\lambda,x)u(\lambda,x)$ satisfies the output LDE 
\begin{equation} \label{eq:OutCC}
\lambda \sigma_2 y(\lambda, x) - \sigma_1 \frac{\partial}{\partial x}y(\lambda,x) +
\gamma_*(x) y(\lambda,x) = 0.
\end{equation}
\end{thm}
\noindent\textbf{Proof:} Let us fix $\lambda\not\in\SPEC(A)$ and a solution $u(\lambda,x)$ of \eqref{eq:InCC}.
Then for $y(\lambda,x) = S(\lambda,x) u(\lambda,x)$ we calculate:
\[ \begin{array}{lll}
\sigma_1 \dfrac{d}{dx} y(\lambda,x)  & = \sigma_1 \dfrac{d}{dx} [(I - C(x)\mathbb X^{-1}(x) (\lambda I - A)^{-1} B(x)\sigma_1) u(\lambda,x)] =\\
& = \sigma_1 \dfrac{d}{dx} u(\lambda,x) - \sigma_1 \dfrac{d}{dx} [ C(x)\mathbb X^{-1}(x) (\lambda I - A)^{-1} B(x)\sigma_1 u(\lambda,x)] \\
& = (\sigma_2\lambda+\gamma)u(\lambda,x) - \sigma_1 \dfrac{d}{dx} [ C(x)\mathbb X^{-1}(x)] \, (\lambda I - A)^{-1} B(x)\sigma_1 u(\lambda,x) \\
& \quad \quad - \sigma_1C(x)\mathbb X^{-1}(x) (\lambda I - A)^{-1} \dfrac{d}{dx} [B(x)] \sigma_1 u(\lambda,x)  \\
&  \quad \quad \quad \quad - \sigma_1C(x)\mathbb X^{-1}(x) (\lambda I - A)^{-1} B(x) \sigma_1 \dfrac{d}{dx} u(\lambda,x) .
\end{array} \]
Using \eqref{eq:DCX-1}, \eqref{eq:DB} and \eqref{eq:InCC} it becomes (notice that $(\lambda I - A)^{-1} B(x)\sigma_1 u(\lambda,x)\in D(A)$)
\[ \begin{array}{llll}
\sigma_1 \dfrac{d}{dx} y(\lambda,x)   =
(\sigma_2\lambda+\gamma)u(\lambda,x) - \\
-  [\sigma_2 C(x)\mathbb X^{-1}(x) A + \gamma_*(x) C(x) \mathbb X^{-1}(x)] (\lambda I - A)^{-1} B(x)\sigma_1 u(\lambda,x) +\\
\quad \quad + \sigma_1C(x)\mathbb X^{-1}(x) (\lambda I - A)^{-1} [A B(x)\sigma_2 + B(x) \cancel{\gamma}] u(\lambda,x) \\
\quad \quad \quad \quad -\sigma_1C(x)\mathbb X^{-1}(x) (\lambda I - A)^{-1} B(x) (\sigma_2\lambda+\cancel\gamma) u(\lambda,x) =
\end{array} \]
Let us combine the last two terms and add $\pm \lambda I$ next to $A$:
\[ \begin{array}{llll}
= (\sigma_2\lambda + \gamma)u(\lambda,x) - \\
-  [\sigma_2 C(x)\mathbb X^{-1}(x) (A \pm \lambda I)+ \gamma_*(x) C(x) \mathbb X^{-1}(x)] (\lambda I - A)^{-1} B(x)\sigma_1 u(\lambda,x) +\\
\quad \quad + \sigma_1C(x)\mathbb X^{-1}(x) (\lambda I - A)^{-1} (A -\lambda I)B(x)\sigma_2 u(\lambda,x) = \\
= (\sigma_2\lambda + \gamma)u(\lambda,x) + \sigma_2 C(x)\mathbb X^{-1}(x)  B(x)\sigma_1 u(\lambda,x) - \\
- [\sigma_2 C(x)\mathbb X^{-1}(x) \lambda + \gamma_*(x) C(x) \mathbb X^{-1}(x)] (\lambda I - A)^{-1} B(x)\sigma_1 u(\lambda,x) - \\
\quad \quad - \sigma_1C(x)\mathbb X^{-1}(x) B(x)\sigma_2 u(\lambda,x) = \\
= (\sigma_2\lambda + \gamma + \sigma_2 C(x)\mathbb X^{-1}(x)  B(x)\sigma_1 -\sigma_1C(x)\mathbb X^{-1}(x) B(x)\sigma_2 )u(\lambda,x) - \\
- [\sigma_2 C(x)\mathbb X^{-1}(x) \lambda + \gamma_*(x) C(x) \mathbb X^{-1}(x)] (\lambda I - A)^{-1} B(x)\sigma_1 u(\lambda,x).
\end{array} \]
Using \eqref{eq:Linkage} and the definition of $S(\lambda,x)$ we obtain that
\[ \begin{array}{llll}
\sigma_1 \dfrac{d}{dx} y(\lambda,x) & = [\sigma_2 \lambda + \gamma_*(x)] u(\lambda,x) - \\
& \quad\quad [\sigma_2 \lambda  - \gamma_*(x)] C(x)\mathbb X^{-1}(x) (\lambda I - A)^{-1} B(x) \sigma_1 u(\lambda,x) = \\
&= [\sigma_2 \lambda + \gamma_*(x)] [I - C(x)\mathbb X^{-1}(x) (\lambda I - A)^{-1} B(x) \sigma_1] u(\lambda,x)=\\
&= (\sigma_2\lambda + \gamma_*(x)) S(\lambda,x) u(\lambda,x) = \\
&= (\sigma_2\lambda + \gamma_*(x)) y(\lambda,x).\qed
\end{array} \]
One of the corollaries of this Theorem is that the function $S(\lambda,x)$ must satisfy \cite{bib:CoddLev} the following differential equation
\begin{equation} \label{eq:DS}
\dfrac{\partial}{\partial x} S(\lambda,x) = \sigma_1^{-1} (\sigma_2\lambda+\gamma_*(x)) S(\lambda,x)- S(\lambda,x) \sigma_1^{-1} (\sigma_2\lambda+\gamma).
\end{equation}
Moreover, defining the fundamental solutions $\Phi(\lambda,x), \Phi_*(\lambda,x)$
\begin{eqnarray}
\label{eq:DefPhi} \lambda \sigma_2 \Phi(\lambda,x) -\sigma_1 \frac{\partial}{\partial x}\Phi(\lambda,x) + \gamma \Phi(\lambda,x) = 0, \quad \Phi(\lambda,0) = I, \\
\label{eq:DefPhi*} \lambda \sigma_2 \Phi_*(\lambda,x) - \sigma_1 \frac{\partial}{\partial x}\Phi_*(\lambda,x) +
\gamma_*(x) \Phi_*(\lambda,x) = 0, \quad \Phi_*(\lambda,0) = I,
\end{eqnarray}
we also obtain that
\begin{equation} \label{eq:SPhi*S0Phi}
S(\lambda,x) = \Phi_*(\lambda,x) S(\lambda,0) \Phi^{-1}(\lambda,x).
\end{equation}

% -------------------------------------------------------------------------------------------------------------
\subsection{\label{sec:preVConcstr}Standard construction of a prevessel}
Now we present the standard construction of a prevessel $\mathfrak {preV}$ from a node $\mathfrak N_0$.
under assumption that the operators $A, A_\zeta$ are generators of analytic semi groups. In general, it is enough to demand that
$A, A_\zeta$ possess ``functional calculus''. We will see later in formula \eqref{eq:PhiForm} that 
\[ \Phi(\lambda,x) = \bbmatrix{\cosh(\lambda x) & i \sinh(\lambda x)\\-i\sinh(\lambda x) & \cosh(\lambda x)} \]
for the Canonical systems parameters.
\begin{defn} Let 
\[ \mathfrak N_0 = \bbmatrix{
C_0 & A_\zeta, \mathbb X_0, A & B_0 & \sigma_1\\
\mathbb C^2 & \mathcal K & \mathbb C^2}
\]
be a node, such that $A, A_\zeta$ and $-A, -A_\zeta$ generate semi-groups (or possess ``functional calculus'') and $D(A)=D(A_\zeta)$. 
The \textbf{standard construction} of the operators $B(x), C(x), \mathbb X(x)$ from the node $\mathfrak N_0$ is as follows
\begin{eqnarray}
\label{eq:StConB} B(x) = \dfrac{1}{2\pi i} \int\limits_\Gamma (\lambda I - A)^{-1} B_0 \Phi^{-1}(\lambda,x-x_0)\sigma_1^{-1} d\lambda, \\
\label{eq:StConC} C(x) = \dfrac{1}{2\pi i} \int\limits_\Gamma \Phi(\lambda,x-x_0) C_0 (\lambda I + A_\zeta)^{-1} d\lambda, \\
\label{eq:StConX} \mathbb X(x) = \mathbb X_0 + \int\limits_{x_0}^x B(y)\sigma_2C(y) dy,
\end{eqnarray}
where the path $\Gamma$ is on the right of the spectrum of the generators, going from $e^{-\theta\infty}$ to $e^{-\theta\infty}$.
\end{defn}
\begin{thm} \label{thm:STConpreVessel}The collection
\[ \mathfrak{preV} = \bbmatrix{C(x) & A_\zeta, \mathbb X(x), A & B(x)&  \sigma_1,\sigma_2,\gamma \\
\mathbb C^2 & \mathcal{K} & \mathbb C^2}
\]
defined by the standard construction from the node $\mathfrak N_0$ is a prevessel, coinciding with $\mathfrak N_0$ for $x=x_0$.
\end{thm}
\Proof The condition $B(x)\sigma_2$ is $A$-regular comes from the definition of $B(x)$. Indeed, for all
$\lambda,x$ $(\lambda I - A)^{-1} B_0 \Phi(\lambda,x-x_0) \in D(A)$. By the existence of the functional calculus, it follows that \eqref{eq:DB},
\eqref{eq:DC} hold. The equation \eqref{eq:DX} is immediate and the Lyapunov equation \eqref{eq:Lyapunov} follows from Lemma \ref{lemma:Redund}.
Finally, we have to show that $\mathbb X(x)(D(A)) \subseteq D(A)$. For each $u\in D(A)$
\[ \mathbb X(x) u = \mathbb X_0 u + \int_{x_0}^x B(y)\sigma_2 C(y) u dy.
\]
Here $\mathbb X_0 u\in D(A)$ by the assumptions on $\mathfrak N_0$. $B(y)\sigma_2\in D(A)$ by the $A$-regularity of $B(y)\sigma_2$. Moreover, since
for each $u\in D(A)$
\[ \dfrac{\partial}{\partial x} B(x)\sigma_1C(x) u = - A B(x)\sigma_2C(x)u - B(x)\sigma_2C(x)A_\zeta u \]
by integrating, we will obtain that
\[ \int_{x_0}^x A B(y)\sigma_2C(y)udy = B_0\sigma_1C_0 u - \dfrac{\partial}{\partial x} B(x)\sigma_1C(x) u + \int_{x_0}^x B(x)\sigma_2C(x)A_\zeta u dy
\]
exists. So, by the closeness of the operator $A$, we obtain that
\[ A \int_{x_0}^x B(y)\sigma_2 C(y) u dy = \int_{x_0}^x A B(y)\sigma_2C(y)udy
\]
exists and $\int_{x_0}^x B(y)\sigma_2 C(y) u dy \in D(A)$. 
\qed

% -------------------------------------------------------------------------------------------------------------
\section{Scattering theory of the operator $L$}
We start from the definition of the vessel parameters, which create an inverse scattering theory of $L$ \eqref{eq:DefL}.
\begin{defn} \label{def:CanSymParam}
The Canonical Systems vessel parameters are defined as follows
\[ \sigma_1 = \bbmatrix{0 & i \\ -i &  0},
   \sigma_2 = I = \bbmatrix{1 & 0 \\ 0 & 1},
   \gamma  = 0 = \bbmatrix{0 & 0 \\ 0 & 0}.
\]
\end{defn}
Expanding the transfer function $S(\lambda,x)$ into a Taylour series around $\lambda=\infty$, we obtain a notion of the moment:
\[ S(\lambda,x) = I - \sum_{n=0}^\infty \dfrac{H_n(x)\sigma_1}{\lambda^{n+1}},
\]
where by the definition \textit{the $n$-th moment} $H_n(x)$ of the function $S(\lambda,x)$ is
\begin{equation} \label{eq:Hn}
H_n(x) = C(x)\mathbb X^{-1}(x)A^nB(x).
\end{equation}
Using the zero moment, for example, we obtain that the so called ``linkage condition'' \eqref{eq:Linkage} is equivalent to
\[ \gamma_*(x) = \gamma + \sigma_2 H_0(x)\sigma_1 - \sigma_1 H_0(x)\sigma_2.
\]
There is also a recurrent relation between the moments $H_n(x)$, arising from \eqref{eq:DS}:
\begin{thm} The following recurrent relation between the moments of the vessel $\mathfrak{V}_{reg}$ holds
\begin{eqnarray}
\label{eq:DHn} 	[\sigma_1^{-1}\sigma_2, H_{n+1}\sigma_1] = (H_n)'_x\sigma_1 - \sigma_1^{-1} \gamma_* H_n\sigma_1.
\end{eqnarray}
\end{thm}
\noindent\textbf{Proof: } Follows from the differential equation \eqref{eq:DS} by plugging $S(\lambda,x) = I - \sum\limits_{n=0}^\infty \dfrac{H_n(x)\sigma_1}{\lambda^{n+1}}$. \qed

Let us investigate more carefully the LDEs \eqref{eq:InCC} and \eqref{eq:OutCC}. Denote $u = \bbmatrix{u_1\\u_2}$, then
\eqref{eq:InCC} becomes
\[ \bbmatrix{\lambda u_1 \\\lambda u_2 } - \bbmatrix{iu_2'\\-iu_1'} = \bbmatrix{0\\0}.
\]
Solving this we obtain that 
\begin{equation} \label{eq:uentries}\left\{ \begin{array}{lll}
u_2' = -i\lambda u_1, \\
u_1'' =i \lambda u_2' = \lambda^2 u_1.
\end{array}\right.
\end{equation}
We can see that actually this equation is equivalent to a second-order differential equation with the spectral parameter $\lambda$:
\begin{equation} \label{eq:InCCequiv}
u_1'' = \lambda^2 u_1.
\end{equation}
In order to analyze \eqref{eq:OutCC}, we denote first moment $H_0(x)= [\pi_{ij}] = \bbmatrix{ a_0 & b_0 \\ c_0 & d_0}$, and as a result, the linkage condition \eqref{eq:Linkage} becomes
\[ \gamma_* = \gamma + \bbmatrix{-i(b_0+c_0) & i(a_0-d_0) \\
 i(a_0-d_0) & i(b_0+c_0)}.
\]
So, if we denote $p=-b_0-c_0$, $q=a_0-d_0$ we will obtain the usual canonical systems equation \eqref{eq:CanSys}:
\[
\lambda y - \bbmatrix{0&i\\-i&0} y + \bbmatrix{i p & iq \\iq &-ip} y = 0.
\]
The term ``scattering matrix'' attached to $S(\lambda,0)$ is understood as follows.
\[ y(\lambda,x)= \bbmatrix{y_1(\lambda,x)\\y_2(\lambda,x)}
\]
of solutions of \eqref{eq:CanSys}, can be derived from \eqref{eq:OutCC} in the following form
\[ y(\lambda,x) = S(\lambda,x) \Phi(\lambda,x) \bbmatrix{1\\0} \left(= \Phi_*(\lambda,x) S(\lambda,0)\bbmatrix{1\\0}\right),
\]
where the fundamental matrix $\Phi(\lambda,x)$, solving \eqref{eq:InCC} with the initial condition $\Phi(\lambda,0)=I$ (-the identity matrix) is
\begin{eqnarray}
\label{eq:PhiForm}\Phi(\lambda,x) = \bbmatrix{\cosh(\lambda x) & i \sinh(\lambda x)\\-i\sinh(\lambda x) & \cosh(\lambda x)},
\end{eqnarray}
where $\cosh(t)=\dfrac{e^{t}+e^{-t}}{2}$, $\sinh(t)=\dfrac{e^{t}-e^{-t}}{2}$. It is a metter of simple calculations to verify that $\Phi(\lambda,x)$
satisfies \eqref{eq:InCC}. This matrix is analytic in $\lambda$ and its structure is known, 
so we can study solutions of \eqref{eq:DefL} or equivalently of \eqref{eq:OutCC}, creating in this manner the (inverse) scattering of $L$ \eqref{eq:DefL}.

Let us examine the recurrence relation \eqref{eq:DHn}. We will research for the simplicity of the
presentation the structure of the first moment $H_0(x)$, but almost the same structure will actually 
apply for all moments. Let us denote
\[ H_1(x) = \bbmatrix{ a_1 & b_1 \\ c_1 & d_1 }.
\]
Then \eqref{eq:DHn} becomes
\[ \begin{array}{ll}
0 = -(H_0)'_x\sigma_1  + \sigma_1^{-1} \gamma_* H_0 \sigma_1 + \sigma_1^{-1}\sigma_2H_1(x)\sigma_1 - H_1(x)\sigma_2 = \\
= \bbmatrix
{ -a_1+d_1 + i(a_0b_0+c_0d_0+b_0')& -b_1-c_1-i(a_0^2+c_0(c_0+b_0)-a_0d_0+a_0')\\
-b_1-c_1 +i (b_0(b_0+c_0)-a_0d_0+d_0^2+d_0') & a_1-d_1-i(a_0b_0+c_0d_0+c_0')} ,
\end{array} \]
From where we obtain:
\[ \left\{ \begin{array}{llllll}
a_1-d_1 = i(a_0b_0+c_0d_0+b_0') = i(a_0b_0+c_0d_0+c_0') \\
b_1+c_1 = i (b_0(b_0+c_0)-a_0d_0+d_0^2+d_0') = -i(a_0^2+c_0(c_0+b_0)-a_0d_0+a_0').
\end{array} \right. \]
Solving the second equalities, envolving the elements of $H_0$, we reach
\[ \begin{array}{llllll}
c_0'=b_0' \\
b_0(b_0+c_0)-a_0d_0+d_0^2+d_0' = -(a_0^2+c_0(c_0+b_0)-a_0d_0+a_0')
\end{array} \]
Since we denoted earlier $p = -c_0-b_0$, $q=a_0-d_0$ together with the assumption $c_0=b_0$ the second equation becomes:
\[ p^2+q^2 + d_0'+a_0'=0.
\]
Since $\dfrac{\tau'}{\tau}=\TR(\sigma_2H_0) = a_0+d_0$ from the last equation we obtain that
\begin{equation} \label{eq:TauForm}
 \dfrac{\tau'}{\tau} = - \int_0^x p^2(y) + q^2(y) dy.
\end{equation}
Since additionally $a_0-d_0=q$ the entries of $H_0$ are as follows
\[ H_0 = \bbmatrix{\dfrac{\dfrac{\tau'}{\tau}+q}{2}&-\dfrac{p}{2} \\ -\dfrac{p}{2} & \dfrac{\dfrac{\tau'}{\tau}-q}{2}}
\]
In the general case similar calculations for \eqref{eq:DHn} with $H_n=\bbmatrix{a_n&b_n\\c_n&d_n}$ results in
\[ \left\{ \begin{array}{lll}
a_{n+1}-d_{n+1} = i(-d_n p + b_n q + b_n') = -i(a_n p + c_n q - c_n'), \\
b_{n+1}+c_{n+1} = -i (-c_n p + a_n q + a_n') = -i(b_np+d_nq-d_n').
\end{array} \right. \]
Solving these equations we find that the moment $H_{n+1}$ is expressed via the previous moment using the following
formulas
\begin{equation} \label{eq:Hn+1Hn} \left\{ \begin{array}{lll}
a_{n+1}-d_{n+1} & = i(-d_n p + b_n q + b_n'), \\
b_{n+1}+c_{n+1} & = -i (-c_n p + a_n q + a_n'), \\
a_{n+1}'+d_{n+1}' & = p(b_{n+1}+c_{n+1}) + q (d_{n+1}-a_{n+1}), \\
b_{n+1}'-c_{n+1}' & = p(d_{n+1}-a_{n+1})-q(b_{n+1}+c_{n+1}).
\end{array} \right. \end{equation}

% ---------------------------------------------------------------------------------------------
\section{\label{sec:Uniqeness}``Uniqueness'' of the scattering data}
Let us consider now the uniqueness of the scattering matrix $S(\lambda,0)$. First we prove the following
\begin{lemma} Suppose that we are given a regular vessel
\[ \mathfrak{V}_{reg} = (C(x), A, \mathbb X(x), A_\zeta, B(x); \sigma_1, 
\sigma_2, \gamma, \gamma_*(x);
\mathcal{H},\mathbb C^3;\mathrm I), \]
realizing coefficients $q(x), p(x)$. Let $S(\lambda,x)$ be its transfer function, defined in \eqref{eq:DefS}.
Let $Y(\lambda)$ be an arbitrary $3\times 3$ matrix function, commuting with the fundamental solution 
$\Phi(\lambda,x)$ of \eqref{eq:InCC}. Then $\widetilde S(\lambda,x) = S(\lambda,x)Y(\lambda)$ is the transfer function of a vessel realizing the same
coefficients $q(x), p(x)$.
\end{lemma}
\noindent\textbf{Proof:} By the definition it follows that
\[ S(\lambda,x) = \Phi_*(\lambda,x) S(\lambda,0) \Phi^{-1}(\lambda,x).
\]
So,
\[ \widetilde S(\lambda,x) = S(\lambda,x) Y(\lambda) = \Phi_*(\lambda,x) S(\lambda,0) \Phi^{-1}(\lambda,x)Y(\lambda) = \Phi_*(\lambda,x) S(\lambda,0) Y(\lambda) \Phi^{-1}(\lambda,x)
\]
and realizes the same coefficients $q(x), p(x)$. By the standard construction, there is a vessel $\widetilde{\mathfrak V}$, whose transfer functions is $\widetilde S(\lambda,x)$.
\qed

Let us investigate the structure of a matrix $Y(\lambda)$, commuting with $\Phi(\lambda,x)$.
Using the form \eqref{eq:PhiForm}, it is easy to conclude that a matrix, which commutes with $\Phi(\lambda,x)$ must be of the form
\begin{equation} \label{eq:YForm}
 Y(\lambda) = a(\lambda) I + b(\lambda) \bbmatrix{0&i\\-i&0},
\end{equation}
by considering the coefficients of pure exponents in $\Phi(\lambda,x)Y(\lambda) = Y(\lambda)\Phi(\lambda,x)$.
\begin{thm}[Uniquness of the scattering matrix]\label{thm:Uniqueness} Suppose that $S(\lambda,x)$, $\widetilde S(\lambda,x)$ are the transfer functions of two regular vessels $\mathfrak V_{reg}, \widetilde {\mathfrak V}_{reg}$,
realizing the same potential $\gamma_*(x)$. Then there exists a matrix $Y(\lambda)\in\R$ such that
\[ \widetilde S(\lambda,x) = S(\lambda,x) Y(\lambda).
\]
\end{thm}
\Proof Let us consider the function $S^{-1}(\lambda,x)\widetilde S(\lambda,x)$. By the definition this functions maps solutions of the input LDE \eqref{eq:InCC} to itself:
\begin{multline*}
 S^{-1}(\lambda,x)\widetilde S(\lambda,x) = \left(\Phi_*(\lambda,x) S(\lambda,0) \Phi^{-1}(\lambda,x)\right)^{-1}\Phi_*(\lambda,x) \widetilde S(\lambda,0) \Phi^{-1}(\lambda,x) =\\
= \Phi(\lambda,x)  S^{-1}(\lambda,0)\widetilde S(\lambda,0) \Phi^{-1}(\lambda,x)
\end{multline*}
Plug here, the formula \eqref{eq:PhiForm} and find conditions so that the coefficients of the exponents $e^{\lambda x}, e^{-\lambda x}$ vanish. This is necessary for
making this function bounded at infinity out of the spectrum of $A$. Then calculations show that actually $S^{-1}(\lambda,0)\widetilde S(\lambda,0)$ must commute with
$\Phi(\lambda,x)$ so that this functions cancels all the singularities at infinity. As a result, by the preceding arguments it must be a function $Y(\lambda)$ of the form
\eqref{eq:YForm}. And we obtain that
\[ S^{-1}(\lambda,0)\widetilde S(\lambda,0) = Y(\lambda),
\]
from where the result follows.
\qed

Another, weaker form of the uniqueness is used later in the text and is presented in the next Lemma. We emphasize that a similar theorem lemma was proved in 
the Sturm-Liouville case in \cite{bib:GenVessel} and in \cite{bib:FaddeyevII} for purely continuous spectrum.
\begin{lemma}\label{lemma:uniquegamma*} Suppose that two functions $S(\lambda,x)$, $\widetilde S(\lambda,x)$ are in class $\I(\sigma_1, \sigma_2,\gamma)$, possessing the same initial value
\[ S(\lambda,0) = \widetilde S(\lambda,0)
\]
and are bounded at a neighborhood of infinity, with a limit value $I$ there.
Then $\widetilde \gamma_*(x) = \gamma_*(x)$.
\end{lemma}
\Proof Suppose that
\[ S(\lambda,x) = \Phi_*(\lambda,x) S(\lambda,0) \Phi^{-1}(\lambda,x), \quad \widetilde S(\lambda,x) = \widetilde \Phi_*(\lambda,x) S(\lambda,0) \Phi^{-1}(\lambda,x).
\]
Then
\[ \widetilde S^{-1}(\lambda,x) S(\lambda,x) = \widetilde \Phi_*(\lambda,x)\Phi^{-1}_*(\lambda,x)
\]
is entire (the singularities appear in $S(\lambda,0) = \widetilde S(\lambda,0)$ only and are cancelled) and equal to $I$ (- the identity matrix) at infinity.
By a Liouville theorem, it is a constant function, namely $I$. So $\widetilde \Phi_*(\lambda,x)\Phi^{-1}_*(\lambda,x) = I$ or
\[ \widetilde \Phi_*(\lambda,x) = \Phi_*(\lambda,x).
\]
If we differentiate, we obtain that $\widetilde \gamma_*(x)=\gamma_*(x)$.
\qed

\begin{thm}[Uniqueness of the moments] \label{thm:momentUnique} Suppose that two sequences of moments $H_n(x)$ and $\widetilde H_n(x)$ are finite, differentiable 
and satisfy \eqref{eq:Hn+1Hn} with analytic $\gamma_*(x)$ and $\widetilde \gamma_*(x)$ respectively. Then from
\[ H_n(0) = \widetilde H_n(0), \quad \forall n=0,1,2,\ldots
\]
it follows that $\gamma_*(x)=\widetilde \gamma_*(x)$. If the infinite system of equations \eqref{eq:Hn+1Hn} has
a unique sequence of solutions $H_n(x)$ for a given $\gamma_*(x)$ and initial values $H_n(0)$ then $H_n(x)=\widetilde H_n(x)$.
\end{thm}
\Proof Let us show by the induction that $H_0^{(n)}(0) = (\widetilde H)_0^{(n)}(0)$ for all $n=0,1,2,\ldots$.
And since these two moments are analytic, the result will follow from the uniqueness of the Taylor series.
For $n=0$, $H_0(0)=\widetilde H_0(0)$ and the basis of the induction follows. Then from \eqref{eq:DHn} it follows that
\[
H_0^{(1)}(x) = \sigma_1^{-1}\sigma_2 H_1 - H_1 \sigma_2\sigma_1^{-1} + \sigma_1^{-1} \gamma_* H_0 - H_0 \gamma\sigma_1^{-1}.
\]
Differentiating again, using \eqref{eq:DHn} for $n=0,1$ and the Linkage condition \eqref{eq:Linkage}, we will obtain that
\begin{multline*}
 H_0^{(1)}(x) = \sigma_1^{-1}\sigma_2 H_1 - H_1 \sigma_2\sigma_1^{-1} + \sigma_1^{-1} \gamma_* H_0 - H_0 \gamma\sigma_1^{-1} =\\
 = \sigma_1^{-1}\sigma_2 H_1 - H_1 \sigma_2\sigma_1^{-1} + \sigma_1^{-1} (\gamma + \sigma_2H_0\sigma_1 - \sigma_1H_0\sigma_2)H_0 - H_0 \gamma\sigma_1^{-1}
 = P_2(H_0(x),H_1(x),H_2(x)) 
\end{multline*}
for a non-commutative polynomial $P_2$ with constant coefficients (depending on $\sigma_1,\sigma_1,\gamma$). This shows that a simple
induction results in
\[ H_0^{(n)}(x) = P_n(H_0(x),H_1(x),\ldots,H_n(x),H_{n+1}(x))
\]
for a non-commutative polynomial $P_n$ with constant matrix-coefficients.
As a result, plugging here $x=0$ and using the condition $H_n(0) = \widetilde H_n(0)$
\begin{multline*}
 H_0^{(n)}(0) = P_n(H_0(0),H_1(0),\ldots,H_n(0),H_{n+1}(0)) = \\
= P_n(\widetilde H_0(0),\widetilde H_1(0),\ldots,\widetilde H_n(0),\widetilde H_{n+1}(0))=
(\widetilde H)_0^{(n)}(0).
\end{multline*}
From here it follows that $H_0(x)=\widetilde H_0(x)$ and hence by the linkage condition \eqref{eq:Linkage} $\gamma_*(x)=\widetilde \gamma_*(x)$. Then
the last statement of the Theorem follows from the uniqueness of solutions. \qed

% -----------------------------------------------------------------------------------------
\section{KdV vessels}
We insert a new variable $t$ to the formulas so that all operators and functions depend now on this variable.
We consider the following notion
\begin{defn} \label{defn:KdVpreV}
The collection of operators and spaces 
\begin{equation} \label{eq:DefKdVpreV}
\mathfrak{preV}_{KdV} = \bbmatrix{C(x,t) & A_\zeta, \mathbb X(x,t), A & B(x,t)&  \sigma_1,\sigma_2,\gamma \\
\mathbb C^2 & \mathcal{K} & \mathbb C^2}
\end{equation}
is called a \textbf{KdV preVessel}, if the following conditions hold:
1. $\mathfrak{preV}_{KdV}$ is a node for all $x,t\in\mathbb R$, 2.
operator $B(x,t)\sigma_2$ is $A^2$-regular, $B(x,t)\gamma$ is $A$-regular 3. $C(x,t), \mathbb X(x,t), B(x,t)$ are differentiable in
both variables, when the other one is fixed, subject to the conditions \eqref{eq:DB}, \eqref{eq:DC}, \eqref{eq:DX} and the following evolutionary equations
(for arbitrary $u\in D(A), v\in D(A)$)
\begin{eqnarray}
\label{eq:DBt} \frac{\partial}{\partial t} B  & = iA \dfrac{\partial}{\partial x} B & = - i A (A B \sigma_2 + B \gamma) \sigma_1^{-1}, \\
\label{eq:DCt} \frac{\partial}{\partial t} C u & = - i \dfrac{\partial}{\partial x} C A_\zeta u & 
= - i \sigma_1^{-1}  (-\sigma_2 C A_\zeta + \gamma C) A_\zeta u, \\
\label{eq:DXt} \frac{\partial}{\partial t} \mathbb X v & =  i (A \dfrac{\partial}{\partial x}\mathbb X - i \dfrac{\partial}{\partial x}\mathbb X A_\zeta+iB\gamma C)v &
= (i A B \sigma_2 C - i B\sigma_2 C A_\zeta+iB\gamma C)v ,
\end{eqnarray}
where $\sigma_2=\sigma_2^*$, $\gamma^*=-\gamma$ are $2\times 2$ matrices. The prevessel $\mathfrak{preV}$ is called symmetric if
$A_\zeta=A^*$ and $C(x,t)=B^*(x,t)$ for all $x,t\in\mathbb R$.
\end{defn}

\begin{defn} \label{def:KdVVessel} The collection of operators, spaces and an open set $\Omega\subseteq\mathbb R^2$ 
\begin{equation} \label{eq:DefKdVVGen}
\mathfrak{V}_{KdV} = \bbmatrix{C(x,t) & A_\zeta, \mathbb X(x,t), A & B(x,t)&  \sigma_1,\sigma_2,\gamma,\gamma_*(x,t) \\
\mathbb C^2 & \mathcal{K} & \mathbb C^2 & \Omega }
\end{equation}
is called a (non-symmetric) \textbf{KdV vessel}, if $\mathfrak{V}_{KdV}$ is a KdV prevessel, $\mathbb X(x,t)$ is invertible on $\Omega$, $\mathfrak{V}_{KdV}$
is also an invertible node.
The $2\times 2$ matrix-function $\gamma_*(x,t)$ satisfies the linkage condition \eqref{eq:Linkage}.
The vessel $\mathfrak{V}_{KdV}$ is called symmetric if $A_\zeta=A^*$ and $C(x,t)=B^*(x,t)$ for all $x,t\in\Omega$.
\end{defn}
\begin{thm} Let $\mathfrak{V}_{KdV}$ be a KdV vessel. Suppose that the moments $H_0,\ldots,H_{n+1}$ are finite and
differentiable, then
\begin{equation}\label{eq:DHntKdV}
\dfrac{\partial}{\partial t} H_n = i \dfrac{\partial}{\partial x} H_{n+1} + i \dfrac{\partial}{\partial x} [H_0] \sigma_1 H_n.
\end{equation}
The transfer function $S(\lambda,x,t)$ \eqref{eq:S0realized} satisfies the following differential equation
\begin{equation} \label{eq:PDEforS}
 \dfrac{\partial}{\partial t} S(\lambda,x,t) = i\lambda \dfrac{\partial}{\partial x} S(\lambda,x,t) +
i \dfrac{\partial}{\partial x} [H_0] \sigma_1 S(\lambda,x,t).
\end{equation}
\end{thm}
\Proof Consider the formula for the moments first.
\[ \begin{array}{lllllll}
\dfrac{\partial}{\partial t} H_n & = \dfrac{\partial}{\partial t}[C\mathbb X^{-1}A^n B] = C_t\mathbb X^{-1}A^nB - 
C \mathbb X^{-1} \mathbb X_t \mathbb X^{-1}A^nB + C\mathbb X^{-1}A^n B_t = \\
& = \text{using evolutionary conditions \eqref{eq:DBt}, \eqref{eq:DXt} } \\
& = C_x(-iA_\zeta) \mathbb X^{-1}A^nB - C \mathbb X^{-1} (iA \mathbb X_x-i\mathbb X_x A_\zeta + iB\gamma C) \mathbb X^{-1}A^nB + C\mathbb X^{-1}A^n (iA) B_x = \\
& = \text{using \eqref{eq:DB}, \eqref{eq:Linkage} and \eqref{eq:Lyapunov}} \\
& = i \dfrac{\partial}{\partial x} H_{n+1} + i \dfrac{\partial}{\partial x} [H_0] \sigma_1 H_n,
\end{array} \]
Similarly one shows the formula \eqref{eq:PDEforS}.
\qed
\begin{cor} The potential $\gamma_*(x,t)$ of a KdV vessel satisfies the following differential equation
\begin{equation} \label{eq:Dgamma*tKdV}
(\gamma_*)_t = - i \gamma_* (H_0)_x\sigma_1 + i \sigma_1 (H_0)_{xx} \sigma_1 +i \sigma_1 (H_0)_x \gamma_*.
\end{equation} 
\end{cor}
\Proof From the linkage condition and \eqref{eq:DHntKdV} for $n=0$ it follows that
\[ \begin{array}{lllllll}
(\gamma_*)_t & = \sigma_2 (H_0)_t \sigma_1 - \sigma_1 (H_0)_t \sigma_2 = \\
& = \sigma_2 [i(H_1)_x  + i (H_0)_x \sigma_1 H_0] \sigma_1 - \sigma_1 [i (H_1)_x + i (H_0)_x \sigma_1 H_0] \sigma_2 \\
& = i \sigma_1 [\sigma_1^{-1} \sigma_2 (H_1)_x - (H_1)_x \sigma_2\sigma_1^{-1}]\sigma_1 + 
i \sigma_2 (H_0)_x \sigma_1 H_0 \sigma_1 - i \sigma_1 (H_0)_x \sigma_1 H_0 \sigma_2.
\end{array} \]
For the first term in this expression we can use the formula \eqref{eq:DHn} for $n=0$, 
then
\[ \begin{array}{lllllll}
(\gamma_*)_t & = i \sigma_1 \dfrac{\partial}{\partial x}[(H_0)_x - \sigma_1^{-1}\gamma_*H_0 + H_0 \gamma \sigma_1^{-1} ]\sigma_1 + 
i \sigma_2 (H_0)_x \sigma_1 H_0 \sigma_1 - i \sigma_1 (H_0)_x \sigma_1 H_0 \sigma_2 \\
& = i \sigma_1(H_0)_{xx}\sigma_1 - i \dfrac{\partial}{\partial x} [\gamma_*H_0 \sigma_1 + \sigma_1 H_0 \gamma]+ 
i \sigma_2 (H_0)_x \sigma_1 H_0 \sigma_1 - i \sigma_1 (H_0)_x \sigma_1 H_0 \sigma_2 \\
& = i \sigma_1(H_0)_{xx}\sigma_1 - i \gamma_* (H_0)_x\sigma_1 - i (\gamma_*)_x H_0 \sigma_1 + \sigma_1 (H_0)_x \gamma +  i \sigma_2 (H_0)_x \sigma_1 H_0 \sigma_1 - i \sigma_1 (H_0)_x \sigma_1 H_0 \sigma_2.
\end{array} \]
Then notice that 
\[ \begin{array}{lllllll}
- i (\gamma_*)_x H_0 \sigma_1 + \sigma_1 (H_0)_x \gamma +  i \sigma_2 (H_0)_x \sigma_1 H_0 \sigma_1 - i \sigma_1 (H_0)_x \sigma_1 H_0 \sigma_2 = \\
= - i [\sigma_2 (H_0)_x\sigma_1 - \sigma_1 (H_0)_x\sigma_2] H_0 \sigma_1 + \sigma_1 (H_0)_x \gamma +  i \sigma_2 (H_0)_x \sigma_1 H_0 \sigma_1 - i \sigma_1 (H_0)_x \sigma_1 H_0 \sigma_2 = \\
= i \sigma_1 (H_0)_x [\gamma + \sigma_2 H_0 \sigma_1 - \sigma_1 H_0 \sigma_2] \\
= i \sigma_1 (H_0)_x \gamma_*,
\end{array} \]
and the result follows. \qed

Let us write explicit formulas for the evolution of $p,q$ arising from \eqref{eq:Dgamma*tKdV}. Since 
\[ \gamma_* = \bbmatrix{ip&iq\\iq&-ip} \]
the (1,1) and (1,2) entries will give some evolutionary formulas. More precisely
\begin{eqnarray}
\label{eq:dpt} \dfrac{\partial}{\partial t} p = q (p^2 + q^2) - \dfrac{q_{xx}}{2}, \\
\label{eq:dqt} \dfrac{\partial}{\partial t} q = - p (q^2 + p^2) + \dfrac{p_{xx}}{2}.
\end{eqnarray}
Multiplying the first expression by $p$, the second by $q$ and summing the results, we obtain after cancellation
\[ pp_t+qq_t = \dfrac{1}{2}(qp_{xx}-pq_{xx}).
\]
Let us denote $\beta(x,t)=\dfrac{\tau'(x,t)}{\tau(x,t)}$. Then from \eqref{eq:TauForm} $\beta'=-(p^2+q^2)$ and the last equality can be rewritten as
\[ \begin{array}{lll}
-\dfrac{\partial}{\partial t} (p^2+q^2) & = qp_{xx}-pq_{xx} \Leftrightarrow \\
-\dfrac{\partial}{\partial t} \beta' & = \dfrac{\partial}{\partial x}[ qp_x-pq_x] \Leftrightarrow \\
-\dfrac{\partial}{\partial x}\dfrac{\partial}{\partial t} \beta & = \dfrac{\partial}{\partial x}[ qp_x-pq_x]  \Leftarrow \\
-\dfrac{\partial}{\partial t} \beta & = qp_x-pq_x,
\end{array} \]
up to a constant of integration. The same formula could be derived as follows: $\beta=a_0+d_0$, and as a result $\beta_t =\TR(H_0)_t$, then 
inserting the formula \eqref{eq:DHntKdV} we will obtain the same result. Differentiating this formula again with respect to $t$ and using
\eqref{eq:dpt}, \eqref{eq:dqt} we obtain
\[ \begin{array}{llll}
-\dfrac{\partial^2}{\partial t^2} \beta & = \dfrac{\partial}{\partial t} [ qp_x-pq_x] =
q_t p_x + q(p_t)_x-p_tq_x -p(q_t)_x \\
& = (- p (q^2 + p^2) + \dfrac{p_{xx}}{2}) p_x + q\dfrac{\partial}{\partial x}\big(q (p^2 + q^2) - \dfrac{q_{xx}}{2}\big)-(q (p^2 + q^2) - \dfrac{q_{xx}}{2})q_x -p\dfrac{\partial}{\partial x}\big(- p (q^2 + p^2) + \dfrac{p_{xx}}{2}\big) \\
& = (p^2+q^2)(p^2+q^2)_x+\dfrac{1}{2}[p''p'-pp'''+q'q''-qq'''],
\end{array} \]
after cancellations. Notice that
\[ - \dfrac{1}{2} \beta'''' = \dfrac{1}{2}(p^2+q^2)''' = pp'''+3p'p''+qq'''+3 q'q''
\]
or
\[ \dfrac{1}{2} \beta'''' + 4 (p''p' + q'q''] = p''p'-pp'''+q'q''-qq'''.
\]
In view of this formula and the definition of $\beta$, we obtain that
\[ -\dfrac{\partial^2}{\partial t^2} \beta = \beta'\beta''+\dfrac{1}{4}\beta''''+\dfrac{\partial}{\partial x} [(p')^2+(q')^2]
\]
or
\begin{equation} \label{eq:Dbetapre}
\dfrac{\partial^2}{\partial t^2} \beta + \beta'\beta''+\dfrac{1}{4}\beta''''+\dfrac{\partial}{\partial x} [(p')^2+(q')^2] = 0.
\end{equation}
In order to express $(p')^2+(q')^2$ in terms of $\beta$ we calculate two expressions:
\[ \begin{array}{llll}
(\dfrac{\partial}{\partial t} \beta)^2 = (qp_x-pq_x)^2 = q^2p_x^2 - 2 pqp_xq_x + p^2q_x^2 \\
\dfrac{1}{4}(\beta'')^2 = (\dfrac{1}{2} (p^2+q^2)')^2 = (pp_x+qq_x)^2 = p^2p_x^2 +2 pqp_xq_x + q^2q_x^2.
\end{array} \]
If we sum them up, we obtain
\begin{equation} \label{eq:px2qx2}
(\dfrac{\partial}{\partial t} \beta)^2 + \dfrac{1}{4}(\beta'')^2 = (p^2+q^2)(p_x^2 + q_x^2) = -\beta' (p_x^2 + q_x^2).
\end{equation}
Finally, plugging \eqref{eq:px2qx2} into \eqref{eq:Dbetapre} we obtain
\[ \dfrac{\partial^2}{\partial t^2} \beta + \beta'\beta''+\dfrac{1}{4}\beta'''' -
\dfrac{\partial}{\partial x} [\dfrac{(\dfrac{\partial}{\partial t} \beta)^2+ \dfrac{1}{4}(\beta'')^2}{\beta'}] = 0,
\]
which is \eqref{eq:CanSysEvol}, after a rearrangement of terms. Thus, we have proved the following Theorem.
\begin{thm} The tau function of a KdV vessel, corresponding to Canonical Systems vessel parameters satisfies equation \eqref{eq:CanSysEvol} on the set
$\Omega$, where the vessel exists.
\end{thm}

% -----------------------------------------------------------------------------
\section{Examples}
\subsection{One dimensional exponential soliton}
A simple soliton is obtained when a one dimensional Hilbert spaces is used for $\mathcal K$. We will use symmetric cases to reduce the number of parameters.
Choose $A=k^2+im$, where $k,m\in\mathbb R$ and define for $b_1=\sqrt{2} k, b_2=k$ the following operators
\[ \begin{array}{llll}
B(x,t) & = \bbmatrix{b_1 e^{A x + i A^2t} + b_2 e^{-A x - i A^2t} & -i b_1 e^{A x + i A^2t} + i b_2 e^{-A x - i A^2t}}, \\
C(x,t) & = B^*(x,t), \\
\mathbb X(x,t) & = -e^{2 k^2 (2 m t - x)} + 2 e^{2 k^2 (-2 m t + x)}.
\end{array} \]
Then
\[ \beta(x,t) = -\dfrac{2 k^2(e^{8 k^2 m t}+2 e^{4 k^2 x}) }{e^{8 k^2 m t}-2 e^{4 k^2 x}}
\]
satisfies \eqref{eq:CanSysEvol}. In this case $\tau(x,t) = \mathbb X(x,t)$. The potential has a ``moving'' singularity, which is
determined from $\tau(x,t)=0$, or $x=2mt-\dfrac{\ln 2}{4k^2}$. Thus $\Omega=\mathbb R^2\backslash\{ (x,t)\mid  x=2mt -\dfrac{\ln 2}{4k^2} \}$.

\subsection{One dimensional rational soliton}
Again, for the one dimensional inner space, suppose that $A=i k$, for $k\in\mathbb R$. Then define
\[ \begin{array}{llll}
B(x,t) & = b \bbmatrix{\cos(kx-k^2t) & \sin(kx-k^2t)}, \\
C(x,t) & = B^*(x,t), \\
\mathbb X(x,t) & = 1 + |b|^2 (x-2kt)
\end{array} \]
Then the collection
\[ \mathfrak{preV} = \bbmatrix{C(x) & A_\zeta, \mathbb X(x), A & B(x)&  \sigma_1,\sigma_2,\gamma \\
\mathbb C^2 & \mathcal{K} & \mathbb C^2}
\]
is a prevessel. On the set $\Omega = \mathbb R^2\backslash\{ (x,t)\mid x=2 k t\}$ it is a vessel, and its beta function
\[ \beta(x,t) = \dfrac{|b|^2}{1+|b|^2 (x-2kt)}
\]
satisfies \eqref{eq:CanSysEvol}. For $t=0$, the vessel has singularity at $x=0$, which ``moves'' with $t$ according to the rule $x=2kt$.

\subsection{Two dimensional soliton}
Let us take $\mathcal K$ to be two dimensional Hilbert space. Define
\[ A = A^* = \bbmatrix{ik_1&0\\0&ik_2}.
\]
Define next
\[ \begin{array}{llll}
B(x,t) & =  \bbmatrix{b_1\cos(k_1x-k_1^2t) & b_1\sin(k_1x-k_1^2t)\\b_2\cos(k_2x-k_2^2t) & b_2\sin(k_2x-k_2^2t)}, \\
C(x,t) & = B^*(x,t), \\
\mathbb X(x,t) & = \bbmatrix{1 + |b_1|^2 (x-2k_1t)& \dfrac{b_1b_2^* \sin[(k_1-k_2)(k_1+k_2)t-x]}{k_2-k_1}\\ 
\dfrac{b_1^*b_2 \sin[(k_1-k_2)(k_1+k_2)t-x]}{k_2-k_1} & 1 + |b_2|^2 (x-2k_2t) }
\end{array} \]
Then
\[ \beta(x,t) = -\dfrac{(k_1-k_2) \left((k_1-k_2) (-|b_2|^2 + |b_1|^2 (-1 + 2 |b_2|^2((k_1+k_2) t-x)))-|b_1b_2|^2 
\sin[2 (k_1-k_2) ((k_1+k_2) t-x)]\right)}{(k_1-k_2)^2 (-1+|b_1|^2 (2 k_1 t-x)) (-1+|b_2|^2(2 k_2 t-x))-|b_1b_2|^2 \sin[(k_1-k_2) ((k_1+k_2) t-x)]^2}
\]
for which it is possible to check that $\beta(x,t) = \dfrac{\tau_x(x,t)}{\tau(x,t)}$, where
\[\tau(x,t)=\det(\mathbb X(x,t)) = (1 - |b_1|^2(2 k_1 t-x)) (1-|b_2|^2(2 k_2 t-x))-\dfrac{|b1 b2|^2 \sin[(k_1-k_2) ((k_1+k_2) t-x)]^2}{(k_1-k_2)^2}.
\]
Again it is possible to check that $\beta(x,t)$ satisfies \eqref{eq:CanSysEvol} by plugging. In this case the set $\Omega$ is defined from $\tau(x,t)$:
\[ \Omega(x,t) = \mathbb R^2\backslash\{ (x,t) \mid \tau(x,t) = 0\}.
\]

\subsection{\label{sec:ClassicScatt}Classical inverse scattering}
In the classical case, it is known that for canonical systems there exists, under condition
\[ \int_{\mathbb R} |p(x)| dx, \int_{\mathbb R} |q(x)| dx < C
\]
there exists a fundamental set of solutions $f_1(\lambda,x)$, $f_2(\lambda,x)$, which satisfy the following asymptotic formulas
for $\Re \lambda = 0$:
\[ f_1(\lambda,x) = e^{\lambda x} \bbmatrix{i\\1} + O(1), \quad
f_2(\lambda,x) = e^{-\lambda x} \bbmatrix{-i\\1} + O(1).
\]
Moreover, we can extend these functions into left half-plane ($\Re \lambda\leq 0$) so that
\begin{equation} \label{eq:f1f2estim}
 f_1(\lambda,x) e^{-\lambda x} = \bbmatrix{i\\1} + O(1),\quad \quad  f_2(\lambda,x) e^{\lambda x} = \bbmatrix{-i\\1} + O(1).
\end{equation}
So, if we define
\[ \Psi_*(\lambda,x) = \bbmatrix{\dfrac{-if_1+if_2}{2} & \dfrac{f_1+f_2}{2} }.
\]
then
\[ \begin{array} {lllll}
\Psi_*(\lambda,x) \Phi^{-1}(\lambda,x) & = \bbmatrix{\dfrac{-if_1+if_2}{2} & \dfrac{f_1+f_2}{2} }
\bbmatrix{\cosh(\lambda x)&i\sinh(\lambda x)\\-i\sinh(\lambda x)&\cosh(\lambda x)}^{-1} \\
& = \bbmatrix{\dfrac{-if_1+if_2}{2} & \dfrac{f_1+f_2}{2} }\bbmatrix{\cosh(\lambda x)&-i\sinh(\lambda x)\\i\sinh(\lambda x)&\cosh(\lambda x)} \\
& = \bbmatrix{\dfrac{-if_1+if_2}{2} & \dfrac{f_1+f_2}{2} }[ \dfrac{e^{\lambda x}}{2} \bbmatrix{1& - i\\i&1} + \dfrac{e^{-\lambda x}}{2} \bbmatrix{1& i\\-i&1}]\\
& = \bbmatrix{-2i \dfrac{f_1 e^{-\lambda x}}{4}+2i\dfrac{f_2 e^{\lambda x}}{4} & 2 \dfrac{f_1 e^{-\lambda x}}{4}+2\dfrac{f_2 e^{\lambda x}}{4}}
\end{array} \]
Finally in view of the estimates \eqref{eq:f1f2estim} for $\Re\lambda\leq 0$ this function behaves for a big $x$ as
\[ S(\lambda,x) = \Psi_*(\lambda,x) \Phi^{-1}(\lambda,x) = \bbmatrix{-\dfrac{i}{2}\bbmatrix{i\\1}+\dfrac{i}{2}\bbmatrix{-i\\1} & 
\dfrac{1}{2} \bbmatrix{i\\1}+\dfrac{1}{2}\bbmatrix{-i\\1}} + O(1) =
I + O(1).
\]
We extend the function $S(\lambda,x)$ to the right half-plane ($\Re\lambda>0$) as follows (using the transpose operation):
\[ S(-\bar\lambda,x) = \dfrac{1}{\det S^*(\lambda,x)}  S^{t*}(\lambda,x), \quad \quad \Re\lambda\leq 0. \]
Then it becomes a globally defined function with a jump along the imaginary axis. It is known that $\det S$ (equal to $a(k)$ for $k=-i\lambda$ in \cite[p. 372]{bib:FaddeyevII}) is actually a nonzero function at the left
half plane and is greater then one on the imaginary axis.

Notice that in this case for each $\lambda\in\mathbb C\backslash i\mathbb R$
it holds
\[ \begin{array} {lllll}
S^*(-\bar\lambda)\sigma_1 S(\lambda) = S^t(\lambda)\sigma_1S(\lambda) =
\dfrac{1}{\det S(\lambda,x)} \Phi^{-1t}(\lambda,x)\Psi^t_*(\lambda,x)  \bbmatrix{0&i\\-i&0} \Psi_*(\lambda,x) \Phi^{-1}(\lambda,x).
\end{array} \]
Then simple calculations show that
\[ \Psi^t_*(\lambda,x)  \bbmatrix{0&i\\-i&0} \Psi_*(\lambda,x) = \bbmatrix{0&i \det\Phi_*\\-i\det\Phi_*&0}.
\]
Notice that $\det S = \det\Phi_* \det\Phi^{-1}=\det\Phi_*$ and we obtain that
\[ \begin{array} {lllll}
S^*(-\bar\lambda)\sigma_1 S(\lambda) = S^t(\lambda)\sigma_1S(\lambda) =
\dfrac{1}{\det S} \Phi^{-1t}(\lambda,x) \bbmatrix{0&i \det S\\-i\det S&0}\Phi^{-1}(\lambda,x) =
\bbmatrix{0&i\\-i&0},
\end{array} \]
using simple calculations for the last equality. Thus we obtain that the function $S(\lambda,x)$ is globally defined
and is $\sigma_1$-symmetric. The initial value $S(\lambda,0)$ possesses the same properties. Let us define
\[ \mathcal H = L^2(\mathbb R)
\]
and $A=i\mu$, the operator of multiplication on $i$ and the variable. This becomes an anti self-adjoint unbounded operator
on $\mathcal H$ with the obvious domain. We represent the function $S(\lambda,0)$ as follows
\[ S(\lambda,0) = I - \int\limits_{\mathbb R} \dfrac{1}{\lambda - i\mu} \bbmatrix{a(\mu)&b(\mu)\\c(\mu)&d(\mu)} \sigma_1 d\mu
\]
Where the matrix $\bbmatrix{a(\mu)&b(\mu)\\c(\mu)&d(\mu)}\sigma_1$ parametrizes the jumps of $S$ along the imaginary axis. It is a very
well known result that a function possessing such jump can be represented in the form above. We define a measure
\[ d\bar\rho = \bbmatrix{a(\mu)&b(\mu)\\c(\mu)&d(\mu)} d\mu
\]
ad defined a space $\mathcal K$, equipped with this measure. Then for $B_0=I=C_0^*$ and $\mathbb X_0=I$ it follows that
\[ S(\lambda,0) = I - C_0 \mathbb X_0^{-1} (\lambda I - A)^{-1}B_0\sigma_1.
\]
applying the standard construction to the obtained node, we will realize the potential $\gamma_*(x)$, from we started
and which satisfied the classical condition.

\subsection{\label{sec:AnalScatt}Scattering of analytic potentials}
Let us start from an analytic potential, in which $p(x),q(x)$ are analytic functions of $x$. Starting from the zero moment
\[ H_0 = \bbmatrix{\dfrac{\dfrac{\tau'}{\tau}+q}{2}&-\dfrac{p}{2} \\ -\dfrac{p}{2} & \dfrac{\dfrac{\tau'}{\tau}-q}{2}}, \quad \quad
\tau=\exp(-\int_0^x [p^2(y)+q^2(y)]dy)
\]
and constructing the moments $H_n(x)$ using formula \eqref{eq:Hn+1Hn} 
\[ \left\{ \begin{array}{lll}
a_{n+1}-d_{n+1} & = i(-d_n p + b_n q + b_n'), \\
b_{n+1}+c_{n+1} & = -i (-c_n p + a_n q + a_n'), \\
a_{n+1}'+d_{n+1}' & = p(b_{n+1}+c_{n+1}) + q (d_{n+1}-a_{n+1}), \\
b_{n+1}'-c_{n+1}' & = p(d_{n+1}-a_{n+1})-q(b_{n+1}+c_{n+1}).
\end{array} \right. \]
we will obtain that their values at zero can be chosen so that $H_n(x)=i^n \bbmatrix{h^n_{11}&h^n_{12}\\h^n_{21}&h^n_{22}}$, where
the functions $h^n_{ij}$ are real valued. This follows immediately from the induction. Define next a $2\times 2$ measure $d\bar\rho$,
satisfying
\[ i^n \int\limits_0^\infty d\bar\rho(\mu) \mu^n = H_n(0) = i^n \bbmatrix{h^n_{11}(0)&h^n_{12}(0)\\h^n_{21}(0)&h^n_{22}(0)}.
\]
the existence of such a measure for each entry follows from \cite{bib:Boas1939}. We define $\mathcal K$ to be the Krein space
of column-functions $f=\bbmatrix{f_1(\mu)\\f_2(\mu)}$, defined by the metric $d\bar\rho(\mu)$ and define
\[ A = i\mu, A_\zeta f = -Af-\sigma_1 \int\limits_0^\infty d\bar\rho(\mu)f(\mu), \quad B_0 = C_0^* = \bbmatrix{1&0\\0&1}, \quad \mathbb X_0 = I.
\]
Obviously, $D(A)=D(A_\zeta)$ and the collection
\[ \mathfrak N = \bbmatrix{
C_0 & A_\zeta, \mathbb X_0, A & B_0 & \sigma_1\\
\mathbb C^2 & \mathcal K & \mathbb C^2}
\]
is an invertible node. Then it is possible to show analogously as it is done in \cite{bib:AnalPotKdV} for the KdV case,
that this node, when evolved realizes the given analytic potential $\gamma_*(x)$. When is further evolved with respect to $t$,
creates by its the logarithmic derivative of the tau function a solution of \eqref{eq:CanSysEvol}. The details are analogous to \cite{bib:AnalPotKdV}.


\begin{thebibliography}{10}



\bibitem{bib:agCanSys}
D.~Alpay and I.~Gohberg.
\newblock Inverse spectral problems for differential operators with rational
  scattering matrix functions.
\newblock {\em J. Differential Equations}, 118:1–--19, 1995.

\bibitem{bib:AGKSCanSys}
D.~Alpay, I.~Gohberg, M.A. Kaashoek, and A.L. Sakhnovich.
\newblock Direct and inverse scattering problem for canonical systems with a
  strictly pseudo – exponential potential.
\newblock {\em Mathematische Nachrichten}, 215(1):5--31, 2000.

\bibitem{bib:amv}
D.~Alpay, A.~Melnikov, and V.~Vinnikov.
\newblock Un algorithme de {S}chur pour les fonctions de transfert des syst\`emes
    surd\'etermin\'es invariants dans une direction
\newblock {\em Comptes-Rendus math\'ematiques (Paris)}, 347(13--14):729--733, 2009.


\bibitem{bib:SchurIEOT}
D.~Alpay, A.~Melnikov, and V.~Vinnikov.
\newblock Schur algorithm in the class {$I$} of {$J$}-contractive functions
  intertwining solutions of linear differential equations.
\newblock {\em IEOT}, 74(3):313--344, 2012.

\bibitem{bib:ArDymCanSys}
D.~Z. Arov and H.~Dym.
\newblock J-inner matrix functions, interpolation and inverse problems for
  canonical systems, {V}: The inverse input scattering problem for {W}iener
  class and rational pxq input scattering matrices.
\newblock {\em Integral Equations and Operator Theory}, 43(1):68--129, 2002.

\bibitem{bib:bgr}
J.~Ball, I.~Gohberg, and L.~Rodman.
\newblock {\em Interpolation of rational matrix functions}.
\newblock Operator Theory: Advances and Applications. Birkh{\" a}user Verlag,
  Basel, 1990.

\bibitem{bib:Boas1939}
Jr. Boas, R.~P.
\newblock The {S}tieltjes moment problem for functions of bounded variation.
\newblock {\em Bull. Amer. Math. Soc.}, 45(6):399–--404, 1939.

\bibitem{bib:Brodskii}
M.S. Brodski{\u i}.
\newblock {\em Triangular and {J}ordan representations of linear operators}.
\newblock Translations of AMS, 1971.

\bibitem{bib:BL}
M.S. Brodskii and M.S. Liv{\v s}ic.
\newblock Spectral analysis of non-self-adjoint operators and intermediate
  systems ({R}ussian).
\newblock {\em Uspehi Mat. Nauk (N.S.)}, 13(1 (79)):3--85, 1958.

\bibitem{bib:CoddLev}
E.~A Coddington and N.~Levinson.
\newblock {\em Theory of ordinary differential equation}.
\newblock Mc-Graw Hill, 1955.

\bibitem{bib:KreinReal}
A.~Diksma, H.~Langer, and H.S.V. de~Snoo.
\newblock Representations of holomorphic operator functions by means of
  resolvents of unitary or self-adjoint operators in {K}rein spaces.
\newblock {\em Operator Theory: Adv. and App.}, 24:123--143, 1987.
\newblock Birkhauser Verlag, Berlin.

\bibitem{bib:EngelNagel}
Klaus-Jochen Engel and Rainer Nagel.
\newblock {\em One-parameter semigroups for linear evolution equations}.
\newblock Springer, 2000.

\bibitem{bib:FaddeyevII}
L.~D. Fadeev.
\newblock The inverse problem in the quantum theory of scattering, {II}
  {English} translation.
\newblock {\em Itogi Nauk. i Techn.}, 4:93--180, 1974.

\bibitem{bib:GGKM}
C.S. Gardner, J.M. Greene, M.D. Kruskal, and R.M. Miura.
\newblock Method for solving the {K}orteweg-de {V}ries equation.
\newblock {\em Phys. Rev. Lett.}, 19:1095--1097, 1967.

\bibitem{bib:GasymobCanSys}
M.~G. Gasymov.
\newblock Inverse problem of scattering theory for a system of {D}irac
  equations of order 2n.
\newblock {\em Tr. Moscow Matem. O-va}, 19:41--112, 1968.

\bibitem{bib:GKSCanSys}
M.~A.~Kaashoek I.~Gohberg and A.~L. Sakhnovich.
\newblock Canonical systems with rational spectral densities: Explicit formulas
  and applications.
\newblock {\em Math. Nachr.}, 149:93–--125, 1998.

\bibitem{bib:KreinCanSys}
M.~G. Krein.
\newblock On integral equations leading to a second order differential equation
  ({R}ussian).
\newblock {\em Dokl. Akad. Nauk SSSR}, 97:21--24, 1954.

\bibitem{bib:KreinMelikCanSys}
M.~G. Krein and F.~E. Melik-Adamyan.
\newblock Theory of {S}-matrices of canonical differential equations with
  summable potential.
\newblock {\em Dokl. Akad. Nauk SSSR}, 46(14):150--155, 1968.

\bibitem{bib:Vortices}
M.S. {L}iv\v sic.
\newblock Vortices of 2{D} systems.
\newblock {\em Operator Theory: Advances and Applications}, 123:7--41, 2001.

\bibitem{bib:GenVessel}
A.~Melnikov.
\newblock On a theory of vessels and the inverse scattering.
\newblock http://arxiv.org/abs/1103.2392.

\bibitem{bib:ENLS}
A.~Melnikov.
\newblock On construction of solutions of the evolutionary {N}on {L}inear
  {S}chr{\" o}dinger equation.
\newblock http://arxiv.org/abs/1209.0179.

\bibitem{bib:KdVVessels}
A.~Melnikov.
\newblock Solution of the {K}d{V} equation using evolutionary vessels.
\newblock http://arxiv.org/abs/1110.3495.

\bibitem{bib:SLVessels}
A.~{M}elnikov.
\newblock Finite dimensional {S}turm {L}iouville vessels and their tau
  functions.
\newblock {\em IEOT}, 71(4):455--490, 2011.

\bibitem{bib:AnalPotKdV}
A.~Melnikov.
\newblock Solution of the {K}d{V} equation on the line with analytic initial
  potential.
\newblock http://arxiv.org/abs/1303.5324, 2013.

\bibitem{bib:MelVinC}
A.~Melnikov and V.~Vinnikov.
\newblock Overdetermined conservative 2{D} systems, invariant in one direction
  and a generalization of {P}otapov's theorem.
\newblock http://arxiv.org/abs/0812.3970.

\bibitem{bib:Staffans}
O.~Staffans.
\newblock {\em Well-Posed linear systems}.
\newblock Encyclopedia of math. Cambridge, 2005.

\end{thebibliography}
\end{document}